\documentclass[11pt]{article}
\usepackage[a4paper,margin=28mm]{geometry}
\usepackage{amsmath,amssymb,amsthm,mathtools}
\usepackage[T1]{fontenc}
\usepackage{lmodern}
\usepackage{microtype}
\usepackage{xcolor}
\usepackage[colorlinks=true,linkcolor=blue!50!black,citecolor=blue!50!black,urlcolor=blue!50!black]{hyperref}
\hypersetup{pdftitle={Projective compactifications of hypersurface complements},pdfauthor={}}
\newtheorem{theorem}{Theorem}[section]
\newtheorem{proposition}[theorem]{Proposition}
\newtheorem{lemma}[theorem]{Lemma}

\theoremstyle{remark}

\newcommand{\C}{\mathbb C}
\newcommand{\Q}{\mathbb Q}
\newcommand{\Z}{\mathbb Z}
\newcommand{\R}{\mathbb R}
\newcommand{\PP}{\mathbb P}
\newcommand{\SH}{SH^{+,S^1}}
\newcommand{\dd}{\mathrm d}
\newcommand{\dc}{\mathrm d^c}
\newcommand{\supp}{\operatorname{supp}}
\newcommand{\Tor}{\operatorname{Tor}}
\newcommand{\Pic}{\operatorname{Pic}}

\newcommand{\sEuler}{\chi_{\mathrm s}}

\title{Projective compactifications of hypersurface complements}
\author{Hairun Wang}
\date{}
\begin{document}
\maketitle
\begin{abstract}
Let $U$ be a smooth affine variety of the form $U=\PP^n\setminus H$, where $H\subset\PP^n$ is a smooth connected hypersurface. We prove that every smooth projective compactification of the complex manifold $U$, by adding a smooth connected divisor, is isomorphic to $\PP^n$. The proof uses the rational Reeb index of the normal circle bundle and positive $S^1$-equivariant symplectic homology. This paper is motivated by the work of Li--Zhou in \cite{li2026kahler} and makes a generalization to their result.\footnotemark{}\footnotetext{The author would like to thank Prof. Zhengyi Zhou and Prof. Chi Li for help discussions.}
\end{abstract}

\section{Introduction}\label{sec:intro}

In general, a noncompact complex manifold $U$ can admit many projective compactifications. Requiring the boundary to be a single smooth divisor imposes a substantial restriction, but it does not in itself explain how the complex structure of the open manifold controls the ambient projective variety. In this paper we prove the following theorem:

\begin{theorem}\label{thm:main}
Let $H\subset\PP^n$ be a smooth connected hypersurface over $\C$. Let $X$ be a smooth projective variety of dimension $n$, and let $Y\subset X$ be a smooth connected divisor. If
\[
X\setminus Y\ \cong\ \PP^n\setminus H
\]
biholomorphically, then $X\cong\PP^n$.
\end{theorem}

When $H$ is a hyperplane, the complement is $\C^n$. In this case, the conclusion establish the standardness of smooth K\"ahler compactifications of $\C^n$, which also follows from the theorem of \cite{li2026kahler}. Their proof relates the normal bundle of the boundary to the indices of periodic Reeb orbits and then uses positive equivariant symplectic homology. Our argument follows their approach for an arbitrary smooth hypersurface complement; the additional issue is that the relevant index may be positive, zero, or negative, and its sign must itself be recovered from the open manifold.\footnotemark{}

We remark that unlike it in \cite[Theorem~1.1]{li2026kahler}, we do not claim the pair isomorphism $(X,Y) \simeq (\PP^n,H)$. Whether this is true remains a separate problem, and as far as to the author's knowledge, it is still open. We refer to \cite{blanc2024complements}, \cite{dubouloz2015log} for partial results.

\footnotetext{\textbf{AI usage}. The author used ChatGPT 6.0. The idea of this paper is orginated from \cite{li2026kahler} and the discussions with Prof. Zhou. The AI contributed to the proposition in Section~\ref{subsec:unbounded} and part of the paper writing. The author have reviewed the final manuscript and take
full responsibility for its mathematical content and proofs.}

\subsection{Outline of the argument}

Write $d=\deg H$ and identify the two complements with a complex manifold $U$. Assume first that $n\geq3$. The Gysin sequence and weak Lefschetz theorem imply that
\[
H^2(X;\Q)=\Q[Y],\qquad Y\text{ is ample}.
\]
Consequently there is a rational number $a$ such that $c_1(TX)=a[Y]$. If $L=N_{Y/X}$, adjunction gives
\[
c_1(TY)=R\,c_1(L),\qquad R=a-1.
\]
For the original compactification, the corresponding coefficient is
\[
R_0=\frac{n+1-d}{d}.
\]
The main step is to prove that $R=R_0$. Once this is established, it follows that $c_1(TX)=\frac{n+1}{d}[Y]$. Integral cohomology then shows that the divisibility of $[Y]$ is exactly $d$, so $X$ has Fano index $n+1$ and is projective space by Kobayashi--Ochiai.

To prove $R=R_0$, we use the fact that the coefficient $R$ has a dynamical interpretation, which appears in the symplectic homology theory and is eventually invariant under deformation. The details will be explained in the later chapters.

Section~\ref{sec:algebra} gives the topological reduction and the integral final step. Section~\ref{sec:stein} constructs the Stein fillings and their regular contact boundaries. Section~\ref{sec:floer} states the Floer-theoretic input, explains its gradings, and constructs the determinant trivialization used in the index calculation. Section~\ref{sec:detection} proves the algebraic support statements. The proof of the main theorem is completed in Section~\ref{sec:finish}.

\subsection{Conventions}

All varieties are over $\C$. A divisor is understood to be reduced, and $[D]$ denotes $c_1(\mathcal O(D))$. Homology and cohomology have rational coefficients unless another coefficient group is displayed. The symbol $\chi$ denotes the ordinary topological Euler characteristic. We use the convention $\dc f=-\dd f\circ J$ for a complex structure $J$; a smooth function is strictly plurisubharmonic when $\dd\dc f$ is positive.

Symplectic homology is taken with rational coefficients and with all free homotopy classes of loops included. We use a homological grading: the differential lowers degree by one. A rational degree and a parity in $\Z/2$ will be used simultaneously. The parity is defined independently of the rational degree; it is not obtained by reducing a rational number modulo two.

\section{Topology and the algebro-geometric reduction}\label{sec:algebra}

We first isolate the information supplied by the topology of the complement. In dimension at least three this reduces the theorem to determining a single rational coefficient. Throughout this section, $U=X\setminus Y\cong\PP^n\setminus H$ and $d=\deg H$.

\begin{lemma}\label{lem:topology}
For $n\geq2$ one has
\[
H_1(U;\Z)\cong\Z/d\Z,\qquad
\Tor H^2(U;\Z)\cong\Z/d\Z,\qquad c_1(TU)=0\ \text{in }H^2(U;\Q).
\]
Moreover $H^1(X;\Q)=0$. If $n\geq3$, then
\[
H^2(X;\Q)=\Q[Y],\qquad \mathcal O_X(Y)\text{ is ample},
\qquad H^1(Y;\Q)=0.
\]
\end{lemma}

\begin{proof}
The homology sequence of $(\PP^n,U)$ and the Thom isomorphism contain
\[
H_2(\PP^n;\Z)\longrightarrow H_0(H;\Z)
\longrightarrow H_1(U;\Z)\longrightarrow H_1(\PP^n;\Z)=0.
\]
The first map is multiplication by $d$, giving the first assertion. The universal coefficient theorem then identifies the torsion subgroup of $H^2(U;\Z)$ with $\operatorname{Ext}(H_1(U;\Z),\Z)=\Z/d\Z$. All these groups are finitely generated: a complement of a smooth divisor retracts onto the complement of a small open tubular neighborhood in the compact manifold.

The section defining $H$ trivializes $\mathcal O_{\PP^n}(d)$ on $U$. Thus the hyperplane class restricts to zero rationally, and the same holds for $c_1(T\PP^n)=(n+1)c_1(\mathcal O(1))$. The inclusion $U\hookrightarrow X$ surjects on $H_1$, by the homology sequence of $(X,U)$ and the Thom isomorphism. Hence $H_1(X;\Q)=0$.

Assume $n\geq3$. Weak Lefschetz gives $H^1(H;\Q)=0$, and the cohomology Gysin sequence
\[
H^0(H;\Q)\xrightarrow{\ [H]\ }H^2(\PP^n;\Q)
\longrightarrow H^2(U;\Q)\longrightarrow H^1(H;\Q)
\]
shows that $H^2(U;\Q)=0$. Applying the same sequence to $(X,Y)$ shows that $H^2(X;\Q)$ is generated by $[Y]$. This class is nonzero, since its intersection with the $(n-1)$st power of an ample class is positive.

If $B$ is an ample divisor on $X$, write $[Y]=t[B]$ in real cohomology. Intersecting with $[B]^{n-1}$ gives $t>0$. All Nakai--Moishezon intersection inequalities for $Y$ therefore follow from those for $B$, so $Y$ is ample. Weak Lefschetz for $Y\subset X$ now gives $H^1(Y;\Q)=H^1(X;\Q)=0$.
\end{proof}

For $n\geq3$, put $L=N_{Y/X}=\mathcal O_X(Y)|_Y$. There is a unique $a\in\Q$ such that
\begin{equation}\label{eq:a}
c_1(TX)=a[Y].
\end{equation}
Adjunction gives
\begin{equation}\label{eq:R}
c_1(TY)=(a-1)c_1(L),\qquad R:=a-1.
\end{equation}
For the original pair $(\PP^n,H)$ the corresponding number is
\begin{equation}\label{eq:R0}
R_0=\frac{n+1-d}{d}.
\end{equation}

The next lemma records the final algebraic step. Its integral assertion is essential: a rational Chern-class identity alone does not determine the Fano index.

\begin{lemma}\label{lem:integral}
Suppose $n\geq2$, $b_2(X)=1$, $Y$ is ample, and
\begin{equation}\label{eq:targetc1}
c_1(TX)=\frac{n+1}{d}[Y]\quad\text{in }H^2(X;\Q).
\end{equation}
Then $X\cong\PP^n$.
\end{lemma}

\begin{proof}
Equation~\eqref{eq:targetc1} and the Nakai--Moishezon criterion imply that $-K_X$ is ample. Thus $X$ is Fano. Kodaira vanishing gives $H^q(X,\mathcal O_X)=0$ for $q>0$, so the exponential sequence gives an isomorphism
\[
c_1:\Pic(X)\xrightarrow{\ \sim\ } H^2(X;\Z).
\]
This group is torsion-free. Indeed, if $T$ is a torsion line bundle, then $T\otimes K_X^{-1}$ is ample. Kodaira vanishing and Riemann--Roch imply
\[
h^0(X,T)=\chi(X,T)=\chi(X,\mathcal O_X)=1.
\]
A nonzero section of $T$ has no zeros: a nonzero effective divisor cannot have zero rational first Chern class, as seen by intersecting with an ample class. Hence $T$ is trivial.

It follows that $\Pic(X)=\Z A$ for a generator $A$ chosen on the ample ray. Write $[Y]=m c_1(A)$ with $m>0$. The integral Gysin sequence gives an injection
\[
H^2(X;\Z)/\Z[Y]\ \hookrightarrow\ H^2(U;\Z)
\]
whose cokernel is a subgroup of $H^1(Y;\Z)$, hence is torsion-free. The source is $\Z/m\Z$, so it is precisely the torsion subgroup of $H^2(U;\Z)$. Lemma~\ref{lem:topology} yields $m=d$.

Since $H^2(X;\Z)$ is torsion-free, \eqref{eq:targetc1} now gives
\[
-K_X\cong A^{\otimes(n+1)}.
\]
The bundle $A$ is ample, since a positive power of it is $\mathcal O_X(Y)$. The Kobayashi--Ochiai characterization of projective space~\cite{kobayashi1973characterizations} therefore implies $X\cong\PP^n$.
\end{proof}

\section{Stein fillings and normal circle bundles}\label{sec:stein}

This section relates a smooth ample divisor to the symplectic geometry of its complement. We first recall the structures involved, then construct an exhaustion with a regular circle bundle at infinity.

\subsection{Liouville domains and contact boundaries}

A \emph{Liouville domain} is a compact manifold $W$ with boundary, equipped with a one-form $\lambda$ such that $\omega=\dd\lambda$ is symplectic and the vector field $Z_\lambda$ defined by
\[
\iota_{Z_\lambda}\omega=\lambda
\]
points outwards along $\partial W$. Its boundary carries the \emph{contact form} $\alpha=\lambda|_{\partial W}$: if $\dim_\R W=2n$, then $\alpha\wedge(\dd\alpha)^{n-1}$ is nowhere zero. The hyperplane bundle $\xi=\ker\alpha$ is symplectic with respect to $\dd\alpha$. The \emph{Reeb vector field} is specified by
\[
\alpha(R_\alpha)=1,\qquad \iota_{R_\alpha}\dd\alpha=0.
\]
The \emph{completion} $\widehat W$ is obtained by adjoining $[1,\infty)\times\partial W$ with one-form $r\alpha$. The Liouville flow identifies a collar of the original boundary with a collar of $\{r=1\}$.

A strictly plurisubharmonic exhaustion $F$ on a complex manifold gives the exact symplectic form $\dd\dc F$. Its regular sublevel sets are Liouville domains with $\lambda=\dc F$. We call the resulting Stein structure \emph{complete} if its Liouville vector field is complete, and \emph{of finite type} if there are no critical points outside a compact set. For such a structure, a sufficiently large sublevel set has completion equal to the whole manifold, up to exact symplectomorphism.

\subsection{The contact boundary of an ample divisor}

\begin{lemma}\label{lem:Stein}
Let $D$ be a smooth ample divisor in a smooth projective manifold $Z$, and put $L=N_{D/Z}$. The complement $V=Z\setminus D$ admits a complete finite-type Stein structure. A sufficiently large regular level of its exhaustion has contact structure isomorphic to that on the unit circle bundle $S(L)$ with curvature a positive representative of $c_1(L)$ and Reeb flow
\[
v\longmapsto e^{-2\pi it}v.
\]
In particular, every simple Reeb orbit has period one for this normalization.
\end{lemma}

\begin{proof}
Choose a positive Hermitian metric on $E=\mathcal O_Z(D)$ and let $s$ be its canonical section. On $V$ set
\[
\psi=-\log\|s\|^2.
\]
This is a proper, bounded-below, strictly plurisubharmonic function. If $D=\{w=0\}$ in a holomorphic coordinate chart and the local frame of $E$ has squared norm $e^{-\varphi}$, then
\begin{equation}\label{eq:logpotential}
\psi=-\log|w|^2+\varphi(z,w).
\end{equation}
Its radial derivative is $-2/|w|+O(1)$, so it has no critical points near $D$. A small compactly supported perturbation makes all critical points nondegenerate while preserving strict plurisubharmonicity. There are then only finitely many critical points.

To obtain completeness, use $F=e^\psi$ and $\lambda=\dc F$. Let $g=\dd\dc\psi(\cdot,J\cdot)$, let $Z_\psi=\nabla_g\psi$, and write $b=\|Z_\psi\|_g^2$. A direct calculation gives
\[
Z_\lambda=\frac{Z_\psi}{1+b},\qquad
0\leq\dd\psi(Z_\lambda)=\frac{b}{1+b}\leq1.
\]
Properness of $\psi$ prevents an integral curve from escaping in finite positive or negative time. Moreover, outside a large sublevel set the vector field is transverse to every level. Every trajectory there meets the chosen level in finite backward time: on each compact region between two regular levels, $\dd\psi(Z_\lambda)$ has a positive minimum. The Liouville flow therefore identifies the exterior with a symplectization, as required.

We describe its contact structure. Equip $L$ with the induced metric and write $M=S(L)$. In a local holomorphic frame, with fiber angle $\theta$, the form
\begin{equation}\label{eq:alpha}
\alpha_0=\frac{1}{4\pi}\bigl(-2\dd\theta+\dc\varphi(z,0)\bigr)
\end{equation}
is well defined on $M$. The change of angle under a change of frame cancels the change in the second term, so these expressions patch globally. Its differential is the pullback of the positive curvature form representing $c_1(L)$. Hence $\alpha_0$ is contact, and its Reeb field is $-2\pi\partial_\theta$.

Choose a smooth tubular identification between a neighborhood of the zero section of $L$ and a neighborhood of $D$ whose differential is complex linear along the zero section. Under fiberwise dilation by a small parameter $\varepsilon$, the pulled-back complex structure converges to the normal-bundle complex structure. In the coordinates of \eqref{eq:logpotential}, the rescaled level $\psi=c$ tends to $S(L)$ as $c\to\infty$, with $\varepsilon=e^{-c/2}$. The pulled-back forms $(4\pi)^{-1}\dc\psi$ on these levels converge in $C^\infty$ to \eqref{eq:alpha}. This follows from the expression $-2\dd\theta+\dc\varphi(z,w)$: after rescaling, its coefficients and all their derivatives tend to the values at $w=0$.

For sufficiently large $c$, the interpolating forms between the pullback and $\alpha_0$ remain contact, since the contact condition is open and $M$ is compact. Gray stability identifies their contact structures. Finally, $\dc F=e^\psi\dc\psi$ on every level of $\psi$, so the same conclusion holds for the complete Stein structure.
\end{proof}

The contact structure in Lemma~\ref{lem:Stein} is the regular Boothby--Wang structure. Its defining feature here is particularly simple: the quotient by the Reeb flow is the divisor $D$, and
\[
\dd\pi:\xi\longrightarrow\pi^*TD
\]
is an isomorphism of symplectic vector bundles. We equip $\xi$ with the complex structure pulled back from $TD$ under this isomorphism.

The actual one-form on a chosen level need not equal $\alpha_0$. This causes no difficulty in using its periodic flow. After the contact identification, write that form as $g\alpha_0$ with $g>0$. In the completed end, whose primitive is $rg\alpha_0$, the graph
\[
r=\frac{C}{g}
\]
has induced form $C\alpha_0$. For $C$ sufficiently large, it bounds a larger Liouville domain with the same completion. Rescaling its primitive by $C^{-1}$ gives a domain with boundary form $\alpha_0$. Such a constant rescaling does not change unfiltered symplectic homology or its gradings.

\subsection{Comparison of the completions}

\begin{proposition}\label{prop:Steincomparison}
For $n\geq3$, the complete Stein structures on $U$ obtained from $(X,Y)$ and $(\PP^n,H)$ are exactly symplectomorphic.
\end{proposition}

\begin{proof}
By Lemma~\ref{lem:topology}, both boundary divisors are ample. Pulling one exhaustion back by the given biholomorphism puts the two structures on the same complex manifold $U$. They are complete and have only finitely many critical points by Lemma~\ref{lem:Stein}.

The complete finite-type Stein comparison theorem states that two such strictly plurisubharmonic exhaustions on a fixed complex manifold determine symplectomorphic symplectic manifolds; see Evans~\cite[Lemma~2.3]{evans2011symplectic} and McLean~\cite[Section~2.1]{mclean2012growth}. Let $f$ be the resulting symplectomorphism, and let $\lambda_0,\lambda_1$ be the primitives. Then $f^*\lambda_1-\lambda_0$ is closed. Since $H^1(U;\R)=0$ by Lemma~\ref{lem:topology}, it is exact.
\end{proof}

The local contact identification in Lemma~\ref{lem:Stein} and the global comparison in Proposition~\ref{prop:Steincomparison} have different roles. The former describes the boundary associated with each divisor in terms of its own normal bundle. The latter compares the completed fillings. We will use this comparison to identify their homology invariants; no contactomorphism between the two boundary circle bundles is needed.

\section{The graded invariant and its regular-circle model}\label{sec:floer}

We now introduce the symplectic invariant used in the proof. The geometric information it extracts is the coefficient $R$ in $c_1(TD)=R\,c_1(L)$. Only a regular circle flow is needed, so its local Floer theory reduces to ordinary homology of the divisor.

\subsection{Positive equivariant symplectic homology}

For a Liouville domain $W$, write
\[
\SH_*(W;\Q)
\]
for its positive $S^1$-equivariant symplectic homology. It is built from Hamiltonian Floer complexes on $\widehat W$. The Hamiltonians are chosen so that their nonconstant periodic orbits near infinity correspond to closed Reeb orbits on $\partial W$. The positive theory quotients out the subcomplex generated by the small interior orbits. The equivariant construction incorporates the circle action that changes the initial point of a parametrized loop. Its differential counts Floer trajectories and lowers the homological degree by one. Taking larger Hamiltonians includes Reeb orbits of arbitrarily large period. We use all free homotopy classes of loops, since the normal fibers need not be contractible in $W$.

We use two established properties of this construction: functorial transfer under exact symplectic embeddings, and the action-filtered Morse--Bott description of periodic orbit families. References for the invariant and transfer maps are Gutt--Hutchings~\cite{gutt2018symplectic} and Hutchings~\cite[Proposition~3.7]{hutchings2024zeta}; the Morse--Bott construction is given by Kwon--van Koert~\cite[Theorem~5.4 and Appendix~B]{kwon2016brieskorn}. The grading conventions and the torsion first-Chern-class case are discussed in Li--Zhou~\cite[Section~2.3 and Remark~3.7]{li2026kahler}.

The Conley--Zehnder index is an integer associated with the path of symplectic linear maps obtained by linearizing a flow along a nondegenerate periodic orbit, after choosing a symplectic trivialization along that orbit. It records the rotation of that path relative to the trivialization. Its extension to paths with degenerate endpoints is the Robbin--Salamon index, normalized to assign index two to one positive full rotation in a complex line. We will need only this index for loops and the standard Morse--Bott shift.

If $c_1(TW)$ is torsion, choose $N>0$ and a nowhere-zero section of $(\det_\C TW)^{\otimes N}$. This determines a rational index: one computes the usual index for the $N$-fold direct sum along an orbit using this determinant trivialization, and divides by $N$. Homotopy classes of trivializations of a trivial complex line bundle form a torsor over $H^1(W;\Z)$. Thus, when $H^1(W;\Q)=0$, the trivialization is unique up to homotopy, and the rational grading is canonical and compatible with replacing $N$ by a multiple. Compatible almost complex structures form a contractible space, so the determinant bundle and its grading do not depend on an auxiliary compatible structure.

There is also the usual parity grading coming from the integer Conley--Zehnder index in an ordinary symplectic trivialization along an orbit. A change of such a trivialization changes that index by an even integer, so parity is well defined. Both gradings are retained:
\[
C=\bigoplus_{k\in\Q,\ \epsilon\in\Z/2}C_k^\epsilon,
\qquad
\partial:C_k^\epsilon\longrightarrow C_{k-1}^{\epsilon+1}.
\]
All these are vector spaces over $\Q$; the second grading does not change the coefficient field.

\begin{proposition}\label{prop:comparison}
For $n\geq3$, the fillings associated with $(X,Y)$ and $(\PP^n,H)$ have isomorphic positive $S^1$-equivariant symplectic homology, preserving both the rational grading and parity.
\end{proposition}

\begin{proof}
Proposition~\ref{prop:Steincomparison} gives an exact symplectomorphism of their completions. We explain why this suffices even though it need not preserve the product coordinates near infinity.

An exact symplectic embedding of Liouville domains induces a transfer map from the homology of the larger domain to that of the smaller one. These maps are functorial, and transfer for adjoining a finite piece of a cylindrical end is an isomorphism after passing to unfiltered homology. Exhaust each completion by radial enlargements of a fixed compact core. The image of the first exhaustion under the symplectomorphism and the second exhaustion are cofinal among compact subsets. Choose alternating inclusions between sufficiently large terms. The induced transfer maps have composites equal to transfer maps within the original exhaustions, which are isomorphisms. Thus the two systems, each consisting of isomorphisms, have isomorphic unfiltered homology.

Each core is a deformation retract of $U$. Lemma~\ref{lem:topology} gives torsion $c_1(TU)$ and $H^1(U;\Z)=0$. The determinant-power trivialization is therefore unique up to homotopy. The symplectomorphism and transfer maps preserve the induced rational grading, as well as the ordinary parity grading. The constant rescalings used to normalize the boundary periods also preserve both.
\end{proof}

\subsection{The determinant trivialization from the normal bundle}

For the remainder of this section, let $D$ be a connected smooth projective manifold of complex dimension $n-1\geq1$, and let $L\to D$ be positive. Assume
\begin{equation}\label{eq:monotone}
c_1(TD)=R\,c_1(L)\quad\text{in }H^2(D;\Q),
\qquad H^1(D;\Q)=0.
\end{equation}
Let $\pi:M=S(L)\to D$ carry the form $\alpha_0$ of \eqref{eq:alpha}, with contact bundle $\xi$. Suppose $(M,\alpha_0)$ is the boundary of a Liouville domain $W$ with torsion $c_1(TW)$ and $H^1(W;\Q)=0$. These conditions hold for the two fillings above when $n\geq3$.

The tautological section of $\pi^*L$ is
\[
\sigma(y,v)=v,\qquad (y,v)\in S(L).
\]
It never vanishes. Choose an integer $N>0$ so that $m=NR$ is integral and so that there is an isomorphism of smooth complex line bundles
\begin{equation}\label{eq:B}
B:L^{\otimes m}\xrightarrow{\ \sim\ }(\det_\C TD)^{\otimes N}.
\end{equation}
Such a choice exists: first clear the denominator in \eqref{eq:monotone}, and then multiply $N$ to kill the torsion difference of the integral first Chern classes. No holomorphic isomorphism is required.

For $x=(y,v)\in M$, let
\[
\Delta_x=\det_\C(\dd\pi|_{\xi_x}):\det_\C\xi_x\longrightarrow\det_\C T_yD.
\]
The \emph{normal determinant trivialization} is the nowhere-zero section
\begin{equation}\label{eq:taunormal}
\tau_{\mathrm{nor}}(y,v)
=(\Delta_x^{\otimes N})^{-1}\bigl(B_y(v^{\otimes m})\bigr)
\quad\text{of }(\det_\C\xi)^{\otimes N}.
\end{equation}
For $m<0$, tensor powers mean powers of the dual vector $v^{-1}\in L_y^*$ characterized by $v^{-1}(v)=1$; for $m=0$ the tensor is $1$. Thus \eqref{eq:taunormal} involves only integral tensor powers.

In local frames $e$ for $L$ and $\eta$ for $\det TD$, write $B(e^{\otimes m})=b\eta^{\otimes N}$ and $v=ze$. If $\widetilde\eta=\Delta^{-1}\eta$, then
\begin{equation}\label{eq:taulocal}
\tau_{\mathrm{nor}}(y,v)=b(y)z^m\widetilde\eta^{\otimes N}.
\end{equation}
This expression will make the winding along a fiber explicit.

\begin{lemma}\label{lem:trivialization}
After replacing $N$ by a common multiple, $\tau_{\mathrm{nor}}$ is homotopic to the determinant-power trivialization induced from $W$.
\end{lemma}

\begin{proof}
The relevant part of the circle-bundle Gysin sequence is
\[
H^1(D;\Q)\longrightarrow H^1(M;\Q)
\longrightarrow H^0(D;\Q)\xrightarrow{\ \cup c_1(L)\ }H^2(D;\Q).
\]
The first group vanishes and the last map is injective, since $D$ is connected and $c_1(L)\ne0$. Hence $H^1(M;\Q)=0$, and therefore $H^1(M;\Z)=0$.

Choose a compatible almost complex structure on $TW|_M$ that preserves $\xi$ and sends the outward Liouville field to $R_{\alpha_0}$. Then
\[
TW|_M\cong\xi\oplus\C,
\]
where the second summand is the radial--Reeb plane, trivialized by the Liouville field. A trivialization of $(\det TW)^{\otimes N}$ restricts, using this splitting, to a section $\tau_{\mathrm{fill}}$ of $(\det\xi)^{\otimes N}$. Enlarge $N$ so that this section and \eqref{eq:taunormal} are both defined. Their ratio is a map $M\to\C^*$. Its homotopy class lies in $H^1(M;\Z)=0$, so the two sections are homotopic.
\end{proof}

\begin{lemma}\label{lem:index}
The rational Robbin--Salamon index of a $p$-fold fiber orbit is $2pR$.
\end{lemma}

\begin{proof}
The Robbin--Salamon index extends the Conley--Zehnder index to symplectic paths with degenerate endpoints. For a unitary loop, it is twice the winding number of its complex determinant. Apply this to the linearized Reeb flow on the $N$-fold direct sum of $\xi$, using the determinant frame \eqref{eq:taunormal}.

Parametrize the orbit by $\gamma_p(t)=(y,e^{-2\pi ipt}v)$, $0\leq t\leq1$. Under $\dd\pi$, the linearized flow on $\xi$ is the identity on the fixed vector space $T_yD$. In that fixed description, \eqref{eq:taulocal} shows that the determinant trivialization itself winds $-pm$ times. Consequently the linearized flow relative to this moving trivialization has determinant winding $pm$. Its index on the $N$-fold direct sum is $2pm$, and division by $N$ gives $2pR$. By Lemma~\ref{lem:trivialization}, this is the grading induced from the filling.
\end{proof}

\subsection{The Morse--Bott model}

The orbits of a fixed period are not isolated: all points of $M$ lie on period-$p$ orbits. This is the Morse--Bott situation, in which a smooth family of critical points replaces an isolated critical point. In the equivariant construction the family is divided by reparametrization, so its quotient is $D$. The stabilizer of a $p$-fold orbit is the finite group $\Z/p$; with rational coefficients it contributes no additional homology. The resulting local group is therefore $H_*(D;\Q)$.

The grading shift is the Robbin--Salamon index minus half the real dimension of the orbit space. Thus a class $x\in H_j(D;\Q)$ at period $p$ has rational degree
\begin{equation}\label{eq:degrees}
|x|=2pR-(n-1)+j,\qquad
\epsilon(x)=j+n-1\pmod2.
\end{equation}
For the parity formula, choose an ordinary symplectic trivialization of $\xi$ along the orbit. The unperturbed return path is a loop and has even Maslov index. The Morse--Bott shift is $-(n-1)$, and a Morse degree $j$ contributes $j$ to the index. This gives the second formula, independently of the rational trivialization.

\begin{proposition}\label{prop:model}
Under the assumptions of \eqref{eq:monotone} and the filling assumptions above, $\SH_*(W;\Q)$ is the homology of a complex with underlying vector space
\begin{equation}\label{eq:model}
C=\bigoplus_{p\geq1} C^{(p)},\qquad
C^{(p)}=\bigoplus_{j=0}^{2n-2}H_j(D;\Q).
\end{equation}
The gradings are given by \eqref{eq:degrees}. The differential lowers the rational degree by one, reverses parity, and strictly lowers $p$.
\end{proposition}

\begin{proof}
Use the action-filtered equivariant Morse--Bott construction of Kwon--van Koert~\cite[Theorem~5.4 and Appendix~B]{kwon2016brieskorn}. The action of a fiber orbit is its period $p$, and Floer differentials do not increase action. At each action level the local equivariant group is $H_*(D;\Q)$ with the shift just computed. The orientation local systems are trivial: the horizontal linearized Reeb flow is complex linear and the radial--Reeb factor has a fixed orientation. The complex-linear criterion in \cite[Lemmas~8.7 and~8.12]{kwon2016brieskorn} therefore applies to every iterate. In the terminology of the equivariant theory, these are good orbit families. The action filtration does not require integral indices or a positive value of $R$. The extension to the torsion first-Chern-class setting, retaining rational indices and parity, is explained in \cite[Remark~3.7]{li2026kahler}.

We spell out the algebraic reduction from local chain groups to local homology, since it is useful when $R=0$. In a chain model for this filtered construction, write the differential as $\partial_0+\partial_<$. Here $\partial_0$ preserves the action level and $\partial_<$ strictly lowers it. Over $\Q$, each local complex splits as its homology plus a contractible complex. Choose the splitting to respect both gradings. The homological perturbation construction transfers the full differential to the sum of the local homologies. Every term of the transferred differential contains a factor $\partial_<$ and hence strictly lowers the positive integer $p$. The series defining the transfer is finite on each element, because each such factor lowers $p$. It therefore defines a differential on the direct sum in \eqref{eq:model}, together with a filtered chain equivalence.

This construction may be performed in the compatible chain models below successive finite action cutoffs and then passed to their increasing union. These are the action-filtered models underlying the Morse--Bott construction; the local reductions introduce no completion in the action direction. Filtered direct limits over a field are exact. Thus \eqref{eq:model} computes the unfiltered group even when infinitely many action levels occur in one rational degree. All operations have the stated bidegrees.
\end{proof}

The usual first-page notation for the same model is
\[
E^1_{p,q}=H_{p+q-2pR+n-1}(D;\Q),\qquad p\geq1.
\]
For the arguments below, it is more convenient to use \eqref{eq:model}. This avoids any separate convergence assertion for an unbounded spectral sequence. We will not need to compute its differential: its degree, parity, and direction in the action filtration are sufficient.

\section{Recovering the adjunction coefficient}\label{sec:detection}

We analyze the graded support of the complex from Proposition~\ref{prop:model}. For a rationally graded vector space $V$, write $\supp V=\{k\in\Q:V_k\ne0\}$. The first lemma detects the magnitude of $R$ once its sign is known. The second determines the direction of unbounded support when the base has nonzero Euler characteristic.

\subsection{Extremal degrees}

\begin{lemma}\label{lem:extrema}
Under the hypotheses of Proposition~\ref{prop:model} the following hold.
\begin{enumerate}
\item If $R>0$, then $\min\supp\SH_*(W)=2R-n+1$.
\item If $R<0$, then $\max\supp\SH_*(W)=2R+n-1$.
\item If $R=0$, the support lies in $[1-n,n-1]$, and $\SH_{1-n}(W)$ is infinite-dimensional.
\item If $R\ne0$, every rationally graded component is finite-dimensional. More strongly, the total dimension in any real interval of length one is bounded by a constant independent of the interval.
\end{enumerate}
\end{lemma}

\begin{proof}
Use the reduced complex of Proposition~\ref{prop:model}. If $R>0$, its least degree is $b=2R-n+1$, represented by $H_0(D)$ in column $p=1$. The differential out of this class vanishes. A possible incoming term must have $p>1$ and
\[
2pR-(n-1)+j=b+1,\quad\text{hence}\quad
j=1-2(p-1)R<1.
\]
If this is a nonnegative integer it is zero, in which case its parity is the same as the target's. This is impossible for a differential. Thus the bottom class survives.

If $R<0$, the largest degree is $c=2R+n-1$, represented by the fundamental class in column $p=1$. No class has degree $c+1$, and no differential from column one can decrease the positive column index. The top class therefore survives.

If $R=0$, formula~\eqref{eq:degrees} gives the asserted bounded support. Its degree $1-n$ part is $\bigoplus_{p\geq1}H_0(D)$. It has no outgoing differential by degree, and its possible incoming group is $\bigoplus_{p\geq1}H_1(D)=0$. It survives in its entirety.

Finally, for $R\ne0$, only a uniformly bounded number of integers $p$ can satisfy
\[
2pR-(n-1)+j\in[t,t+1],\qquad 0\leq j\leq2n-2.
\]
The total dimension of each column is the fixed number $\sum_j b_j(D)$. This proves the uniform bound for the complex, and hence for its homology.
\end{proof}

\subsection{The direction of unbounded support}\label{subsec:unbounded}

The support bounds alone do not distinguish a positive index from a negative one: a graded vector space could be bounded in both directions. The following argument excludes this possibility when $\chi(D)\ne0$. It uses the parity grading, which remains available even when the rational degrees are nonintegral.

For a finite-dimensional $\Z/2$-graded vector space $V$, write
\[
\sEuler(V)=\dim V^{\overline0}-\dim V^{\overline1}.
\]
For a finite-dimensional complex with odd differential, this signed Euler characteristic equals that of its homology. Indeed, boundaries and their preimages occur in equal-dimensional pieces of opposite parity. In our application the underlying spaces remain over $\Q$; the superscripts record only parity.

\begin{lemma}\label{lem:tail}
Suppose $\chi(D)\ne0$. If $R>0$, the support of $\SH_*(W)$ is unbounded above. If $R<0$, it is unbounded below.
\end{lemma}

\begin{proof}
Each full column of \eqref{eq:model} has signed Euler characteristic
\begin{equation}\label{eq:columnchi}
\eta=(-1)^{n-1}\chi(D)\ne0.
\end{equation}
Assume first that $R>0$, and set $C^{\leq T}=\bigoplus_{k\leq T}C_k$. Since the differential lowers degree, this is a finite-dimensional subcomplex. Apart from a uniformly bounded number of columns intersecting the upper cutoff, it consists of full columns. Therefore
\begin{equation}\label{eq:Euler-growth}
\sEuler(C^{\leq T})=\frac{\eta}{2R}\,T+O(1)
\quad\text{as }T\longrightarrow+\infty.
\end{equation}
The bounded error follows because the spread of degrees in each column is $2n-2$.

Suppose, to the contrary, that $H_*(C)$ were supported in degrees at most $B$. It is also bounded below and has uniformly bounded dimension in unit intervals, so its total dimension is finite. For $k\leq T-1$, the natural map
\[
H_k(C^{\leq T})\longrightarrow H_k(C)
\]
is an isomorphism: all relevant chains in degrees $k-1,k,k+1$ have been retained. The only additional homology of $C^{\leq T}$ can occur in $(T-1,T]$, where its dimension is bounded by the corresponding chain groups, uniformly in $T$. It follows that $|\sEuler(H(C^{\leq T}))|$ is bounded independently of $T$. This contradicts \eqref{eq:Euler-growth}, since signed Euler characteristic is unchanged on taking homology of a finite-dimensional complex with odd differential.

For $R<0$, use the finite-dimensional quotient
\[
C^{\geq -T}:=C\big/\bigoplus_{k<-T}C_k.
\]
Its signed Euler characteristic is $\eta T/(2|R|)+O(1)$. Away from the bottom interval $[-T,-T+1)$ its homology agrees with that of $C$. If $H_*(C)$ were bounded below, the same uniform-bound argument would contradict this linear growth.
\end{proof}

\subsection{Comparison with a hypersurface in projective space}

We apply the preceding lemmas first to the known compactification. This is enough to determine the sign for the unknown one; no nonvanishing assumption on $\chi(Y)$ is needed.

\begin{lemma}\label{lem:EulerH}
If $n\geq3$, every smooth degree-$d$ hypersurface $H\subset\PP^n$ has $\chi(H)\ne0$.
\end{lemma}

\begin{proof}
Let $h=c_1(\mathcal O_{\PP^n}(1))$. The tangent-normal sequence gives
\[
c(TH)=\frac{(1+h)^{n+1}}{1+dh}\Big|_H,
\qquad \int_H h^{n-1}=d.
\]
Taking the top Chern number and summing the resulting binomial expression yields
\begin{equation}\label{eq:chiH}
\chi(H)=n+1+\frac{(1-d)^{n+1}-1}{d}.
\end{equation}
If $n$ is odd, this number is positive. If $n$ is even, then $n\geq4$; for $d=1,2$ it equals $n$. For $d\geq3$,
\[
(d-1)^{n+1}>(n+1)d-1.
\]
Indeed the difference is positive at $d=3$, since $2^{n+1}>3n+2$ for $n\geq4$, and its derivative as a real function of $d$ is positive for $d\geq3$. Equation~\eqref{eq:chiH} is consequently negative in this last case.
\end{proof}

\begin{proposition}\label{prop:recover}
For the two compactifications in Theorem~\ref{thm:main}, if $n\geq3$, one has $R=R_0$.
\end{proposition}

\begin{proof}
Apply Lemmas~\ref{lem:extrema} and~\ref{lem:tail} to the filling associated with $H$. Lemma~\ref{lem:EulerH} supplies the nonzero Euler characteristic required in the second lemma. Proposition~\ref{prop:comparison} identifies its graded homology with that of the filling associated with $Y$.

If $R_0>0$, the common homology is bounded below and unbounded above. A negative value of $R$ would bound it above, and $R=0$ would bound it in both directions. Thus $R>0$. Comparing the least supported degrees gives
\[
2R-n+1=2R_0-n+1,
\]
so $R=R_0$.

If $R_0<0$, the common homology is bounded above and unbounded below. The support bounds rule out $R\geq0$, and comparison of the greatest supported degrees again gives $R=R_0$.

If $R_0=0$, the degree $1-n$ group is infinite-dimensional. For every nonzero $R$, all graded components are finite-dimensional. Hence $R=0=R_0$.
\end{proof}

\section{Proof of the main theorem}\label{sec:finish}

\subsection{Dimensions at least three}
By \eqref{eq:R}, \eqref{eq:R0}, and Proposition~\ref{prop:recover},
\[
a-1=\frac{n+1-d}{d},\qquad
c_1(TX)=\frac{n+1}{d}[Y].
\]
Lemma~\ref{lem:topology} supplies $b_2(X)=1$ and ampleness of $Y$. Lemma~\ref{lem:integral} now gives $X\cong\PP^n$.

\subsection{Surfaces}
Here the topology of the end determines the numerical information needed for adjunction. This argument includes a cubic boundary, whose Euler characteristic vanishes.

Let $M_Y$ and $M_H$ be the boundaries of small tubular neighborhoods of the two curves. A complement of a smooth connected divisor has one collared end, with cofinal neighborhoods diffeomorphic to $M_D\times(0,\varepsilon)$. The biholomorphism of complements is a proper homeomorphism and identifies their ends. Consequently
\begin{equation}\label{eq:endH1}
H_1(M_Y;\Z)\cong H_1(M_H;\Z).
\end{equation}
One can see this without identifying the cross-sections themselves: the inverse system of first homology groups of a collared end is constant, since the collar inclusions are homotopy equivalences. The two cofinal systems define isomorphic constant systems, hence have isomorphic groups.

For an oriented circle bundle over a genus-$g$ surface with Euler number $e$, its first homology is
\begin{equation}\label{eq:circleH1}
H_1(M;\Z)\cong
\begin{cases}
\Z^{2g}\oplus\Z/|e|\Z,& e\ne0,\\
\Z^{2g+1},& e=0.
\end{cases}
\end{equation}
For example, this follows by abelianizing the usual presentation with a central fiber generator $t$ and relation $\prod_{i=1}^g[a_i,b_i]=t^e$.

For $H\subset\PP^2$, the Euler number is $H^2=d^2$, and its genus satisfies
\[
2g(H)-2=d(d-3).
\]
The rank of $H_1(M_H;\Z)$ is even. Equations~\eqref{eq:endH1}--\eqref{eq:circleH1} therefore rule out $Y^2=0$, and give
\begin{equation}\label{eq:surfaceend}
g(Y)=g(H),\qquad |Y^2|=d^2.
\end{equation}
By Lemma~\ref{lem:topology}, $b_1(X)=0$, and Poincar\'e duality gives $b_3(X)=0$. Additivity of compactly supported Euler characteristic gives
\[
\chi(X)-\chi(Y)=\chi(U)=3-\chi(H).
\]
Here ordinary and compactly supported Euler characteristics of the smooth real four-manifold $U$ agree by Poincar\'e duality. Since the two curves have the same genus, $\chi(X)=3$, whence $b_2(X)=1$.

An effective nonzero divisor on a projective surface with $b_2=1$ has cohomology class a positive multiple of an ample class, and is ample by Nakai--Moishezon. In particular $Y^2>0$, so \eqref{eq:surfaceend} gives $Y^2=d^2$. Write $c_1(TX)=a[Y]$. Adjunction then gives
\[
c_1(TX)\cdot Y=Y^2+2-2g(Y)=3d,
\qquad a=\frac3d.
\]
Thus \eqref{eq:targetc1} holds for $n=2$, and Lemma~\ref{lem:integral} gives $X\cong\PP^2$.

\subsection{Curves}
If $n=1$, a smooth connected hypersurface in $\PP^1$ is a single point, so $U\cong\C$. The divisor $Y$ is also one point. A compact Riemann surface punctured once has first homology of rank twice its genus. Since $H_1(\C;\Z)=0$, the genus of $X$ is zero and $X\cong\PP^1$. This completes the proof of Theorem~\ref{thm:main}.\hfill$\square$

\bibliographystyle{alpha}
\bibliography{refs} 

@article{blanc2024complements,
  title={Complements of hypersurfaces in projective spaces},
  author={Blanc, J{\'e}r{\'e}my and Poloni, Pierre-Marie and Van Santen, Immanuel},
  journal={Journal de l’{\'E}cole polytechnique—Math{\'e}matiques},
  volume={11},
  pages={733--768},
  year={2024}
}

@article{dubouloz2015log,
  title={Log-uniruled affine varieties without cylinder-like open subsets},
  author={Dubouloz, Adrien and Kishimoto, Takashi},
  journal={Bulletin de la Soci{\'e}t{\'e} Math{\'e}matique de France},
  volume={143},
  number={2},
  pages={383--401},
  year={2015}
}

@article{hutchings2024zeta,
  title={Zeta functions of dynamically tame Liouville domains},
  author={Hutchings, Michael},
  journal={arXiv preprint arXiv:2402.07003},
  year={2024}
}

@article{evans2011symplectic,
  title={Symplectic mapping class groups of some Stein and rational surfaces},
  journal={},
  author={Evans, Jonathan David},
  year={2011}
}

@article{mclean2012growth,
  title={The growth rate of symplectic homology and affine varieties},
  author={McLean, Mark},
  journal={Geometric And Functional Analysis},
  volume={22},
  number={2},
  pages={369--442},
  year={2012},
  publisher={Springer}
}

@article{kwon2016brieskorn,
  title={Brieskorn manifolds in contact topology},
  author={Kwon, Myeonggi and Van Koert, Otto},
  journal={Bulletin of the London Mathematical Society},
  volume={48},
  number={2},
  pages={173--241},
  year={2016},
  publisher={Oxford University Press}
}

@article{li2026kahler,
  title={K{\"a}hler compactification of $\mathbb{C}^n$ and Reeb dynamics},
  author={Li, Chi and Zhou, Zhengyi},
  journal={Inventiones mathematicae},
  volume={243},
  number={1},
  pages={181--208},
  year={2026},
  publisher={Springer}
}

@article{kobayashi1973characterizations,
  title={Characterizations of complex projective spaces and hyperquadrics},
  author={Kobayashi, Shoshichi and Ochiai, Takushiro},
  journal={Journal of Mathematics of Kyoto University},
  volume={13},
  number={1},
  pages={31--47},
  year={1973},
  publisher={Duke University Press}
}

@article{gutt2018symplectic,
  title={Symplectic capacities from positive S1--equivariant symplectic homology},
  author={Gutt, Jean and Hutchings, Michael},
  journal={Algebraic \& geometric topology},
  volume={18},
  number={6},
  pages={3537--3600},
  year={2018},
  publisher={Mathematical Sciences Publishers}
}

\end{document}